\documentclass[12pt,reqno]{amsart}
\usepackage{amsmath, amsthm, amssymb, color, mathtools}
\usepackage{enumerate}
\usepackage{comment}
\usepackage[all]{xy}
\usepackage[new]{old-arrows}

\usepackage[T1]{fontenc} 
\usepackage[utf8]{inputenc}

\DeclareUnicodeCharacter{0261}{\k{a}}
\DeclareUnicodeCharacter{0260}{\k{A}}
\DeclareUnicodeCharacter{0281}{\k{e}}
\DeclareUnicodeCharacter{0280}{\k{E}}

\newtheorem{theorem}{Theorem}[section]
\newtheorem{proposition}[theorem]{Proposition}
\newtheorem{lemma}[theorem]{Lemma}

\newtheorem{corollary}[theorem]{Corollary}

\theoremstyle{definition}
\newtheorem{definition}[theorem]{Definition}

\newtheorem{example}[theorem]{Example}
\newtheorem{remark}[theorem]{Remark}

\newcommand{\fra}[1]{{\color{red}{#1}}}

\newcommand{\skipit}[1]{{}}
\newcommand{\prfend}{\hbox to7pt{\hfil}
\par\vskip-\baselineskip\hbox to\hsize
{\hfil\vbox {\hrule width6pt height6pt}}\vskip\baselineskip}

\newcommand{\CC}{\mathbb{K}}

\newcommand {\PP}{\mathbb{P}}
\newcommand {\p}{\mathbb{P}}

\newcommand{\psr}{partially symmetric rank }
\newcommand{\cL}{\mathcal{L}}

\newcommand{\ot}{\otimes}

\newcommand{\ru}[1]{\left\lceil #1\right\rceil}

\DeclareMathOperator{\rk}{srk}
\DeclareMathOperator{\trk}{rk}

\DeclareMathOperator{\Ann}{Ann}
\DeclareMathOperator{\ann}{Ann}
\DeclareMathOperator{\GL}{GL}

\DeclareMathOperator{\Var}{V}
\DeclareMathOperator{\prk}{psrk}

\DeclareMathOperator{\Sym}{Sym}

\newcommand{\Span}[1]{\langle#1\rangle}
\newcommand{\myarrow}[2]{\hbox to #1pt{\hfil$\to$\hfil}{\hskip-#1pt{\raise
10pt\hbox to#1pt{\hfil$\scriptscriptstyle #2$\hfil}}}}

\title{Symmetric rank of some reducible forms}

 \author[L. Colarte-G\'omez]{Liena Colarte-G\'omez} 
 \address{Institute of Mathematics of the Polish Academy of Sciences, \'Sniadeckich 8, 00-656 Warwaw, Poland} \email{lcolartegomez@impan.pl (ORCID 0000-0002-8767-0327)}

 \author[F. Galuppi]{Francesco Galuppi}
 \address{Faculty of Mathematics, Informatics, and Mechanics, University of Warsaw, Banacha 2, 02-097 Warsaw, Poland}
\email{galuppi@mimuw.edu.pl (ORCID 0000-0001-5630-5389)}
\thanks{The second author is the corresponding author.}

\begin{document}
\begin{abstract} In this paper, we study the symmetric rank of products of linear forms and an irreducible quadratic form
. The main result presents a new, non-trivial lower bound for the rank, and the arguments rely on the apolarity lemma. 
In the special case of degree 4, we give a 
list of normal forms for such quartics, and we apply our general result to compute the rank of almost all of them. 
These families of quartics provide examples of polynomials of generic, supergeneric and even maximal rank, as well as an unexpected number of decompositions.
\end{abstract}



\maketitle


\section{Introduction}

Tensors are central objects in pure and applied mathematics, ubiquitous in many branches of science. Simply speaking, a tensor is a multidimensional array, for instance, a way to store large amounts of data. 
A complete list of applications of tensors would be far too long, so we refer the reader to \cite[Section 1.3]{Lan}. Compelled to understand the structure and essential properties of the tensor, we would attempt to decompose it as a combination of simpler tensors.

The \emph{rank} of a tensor $T\in V_1\ot\dots\ot V_s$ is a natural number that indicates the minimal complexity of such a decomposition. For $s=2$, it is the same as the matrix rank. However, tensor rank becomes rather complicated for $s\ge 3$, and the literature contains several variants of rank, tailored for different purposes. Examples include the border rank \cite[Section 1.2.6]{Lan}, the multilinear rank \cite[Definition 2.3.1.4]{Lan} and the geometric rank \cite{geometricrank}. 

This paper focuses on the \emph{symmetric} and the \emph{partially symmetric rank}. Given a field $\CC$ of characteristic 0, we regard the space $\Sym^{d}(\CC^{n+1})$ of  symmetric tensors of order $d$ as the space $\CC[x_0,\hdots,x_n]_d$ of degree $d$ homogeneous polynomial. 
Determining the symmetric rank is a long-standing challenging problem, also called {\em Waring problem} in the literature. The solution is known for quadratic forms, binary forms
, monomials
, and generic forms
. However, computing the rank of a given polynomial is a difficult task: according to \cite[Theorem 8.2 and Conjecture 13.2]{nphard}, it is an NP-hard problem. The comparison between the symmetric and the \psr of a symmetric tensor is also motivated by Comon's conjecture \cite{ComonConj} and its partially symmetric analog \cite{PartComon}. Despite many attempts to prove or disprove it, this interesting conjecture remains open \cite{comonisopen}
.

In \cite{reduciblecubics}, the authors study the symmetric rank of reducible cubics, namely polynomials that are the product of a linear form $\ell_0$ and a quadratic form $q$. Their result inspired us to study the rank of products of linear forms and an irreducible quadratic form, that is, polynomials of the form $\ell_0\cdots \ell_k \cdot q$. We chose this family because it contains examples of polynomials of generic, supergeneric, and even maximal rank that are typically difficult to find otherwise. When the field $\CC$ is algebraically closed of characteristic 0, the main tool to compute the symmetric rank of a polynomial is the celebrated Apolarity lemma, proven in \cite{iarrobinoKanev99}. While it may be sometimes difficult to use in practice, the Apolarity lemma can be applied nicely to our family of polynomials, especially in low degrees.



Our primary purpose is to provide a new, non-trivial lower bound for the symmetric and partially symmetric ranks. For our family of polynomials $\ell_0\cdots \ell_k \cdot q$, we provide such bound in Theorem \ref{Proposition: generalized lower bound} and Corollary \ref{Corollary: geometric conse generalized lower bound}, under a mild assumption. We illustrate the result through examples and compare it with known bounds (Examples \ref{Example: non asymtotic in n} and \ref{Example: high degrees examples}).

In Section \ref{Section: The symmetric rank of quartics} we focus on the symmetric rank of quartics $F = \ell_0\ell_1 \cdot q$. Up to $\GL$-equivalence, this class of polynomial admits a list of ten normal forms, that we compute in Proposition \ref{Prop: normal forms of reducible quartics}. In Proposition \ref{Theorem: lower rank of quartics}, we apply our general lower bound to show that if $\ell_0$ and $\ell_1$ are not proportional, then the symmetric rank of $\ell_0\ell_1q$ is at least $4(n-1)$ and the \psr is at least $3(n-1)$. Finally, we prove an upper bound on the rank of such forms: In Theorem \ref{Prop: upper bounds for reducible quartics}, we establish that the symmetric rank is at most $4(n-1)+2$. For $n=3$ we improve this bound by one. We end these notes by posting some problems and questions.

\vspace{0.5cm}
\noindent{\bf Acknowledgements.} This project started within the framework of the algebraic geometry working group {\em Secant varieties} at IMPAN. The authors would like to warmly thank Jaros\l{}aw Buczy\'{n}ski, Maciej Gałązka, Iryna Raievska, and Maryna Raievska for very useful discussions and comments.

Galuppi acknowledges support by the National Science Center, Poland, under the projects ``Complex contact manifolds and geometry of secants'', 2017/26/E/ST1/00231, and ``Tensor rank and its applications to signature tensors of paths'', 2023/51/D/ST1/02363.

\section{Ranks and apolarity}\label{section: preliminaries}

This preliminary section 
gathers some basic notions and results we shall use in the paper. Let $\CC$ be an algebraically closed field of characteristic zero, and let $R = \CC[x_0,\hdots, x_n]$ be the polynomial ring with the standard grading. 
Given a subset $E \subset R$, we write $\Var(E)$ for the zero-locus $\{p \in \PP^n \mid f(p) = 0\mbox{ for every } f \in E\}$ of $E$. 

\begin{definition}\label{definition: sym rank} Let $F\in R$ be a form of degree $d$. A \emph{symmetric decomposition} (also called  \emph{Waring decomposition}) of $F$ of length $r$ is a set of $r$ linear forms $\{\ell_1,\dots,\ell_r\}\subset R_1$ such that
\[
F\in\langle \ell_1^d,\dots,\ell_r^d\rangle,
\]
that is, $F$ can be decomposed as a linear combination of $\ell_1^{d},\dots,\ell_r^{d}$. The \emph{symmetric rank} (or \emph{Waring rank}) of $F$ is the minimal length of a symmetric decomposition of $F$. We denote it by $\rk(F)$.
\end{definition}


An important tool for computing the symmetric rank of a given polynomial is apolarity. Consider another polynomial ring $S=\CC[y_0,\dots,y_n]$ and let $S$ act on $R$ by differentiation: if $F\in R$ then 
\[
y_i[F] =\frac{\partial F}{\partial x_i}.
\]
The \emph{annihilator} of $F$, also called the \emph{inverse system} or the \emph{apolar ideal} of $F$, is the ideal
\[
\Ann(F)=\{\delta\in S\mid \delta[F]=0\}\subset S
\]
of all differential polynomials which annihilate $F$. The connection to the symmetric rank is the celebrated apolarity lemma, proven in \cite[Theorem 5.3]{iarrobinoKanev99}.

\begin{lemma}[Apolarity lemma]\label{lemma: apolarity} Let $F\in R_d$ be a form of degree $d$, let $[\ell_1],\dots,[\ell_r]\in \PP(R_1)$ be $r$ distinct linear forms and $Z=\{\ell_1^{\vee},\dots,\ell_r^{\vee}\}\subset \PP^n$. Then $F$ can be decomposed as a linear combination of $\ell_1^d,\dots,\ell_r^d$ if and only if $I(Z)\subset \Ann(F)$.

In particular, $\rk(F)\le r$ if and only if $\Ann(F)$ contains the ideal of a set of $r$ points. 
\end{lemma}

In this paper, we are also interested in the partially symmetric rank. 

\begin{definition}\label{definition: partially sym rank} Let $F\in R_d$ be a form of degree $d$. A \emph{partially symmetric decomposition} of $F$ of length $r$ is a simultaneous decomposition of all partial derivatives of $F$, that is, a set of $r$ linear forms $\{\ell_1,\dots,\ell_r\}\subset R_1$ such that
\[
y_i[F]\in\langle \ell_1^{d-1},\dots,\ell_r^{d-1}\rangle\mbox{ for every } i\in\{0,\dots,n\}.
\]
The \emph{partially symmetric rank} of $F$ is the minimal length of a partially symmetric decomposition of $F$. We denoted it by $\prk(F)$.
\end{definition}

\begin{remark}
\label{remark: sym rank bigger than par sym rank}    Every symmetric decomposition of $F$ is a partially symmetric decomposition of $F$. Thus, $\prk(F) \leq \rk(F)$ for any form $F\in R$.
\end{remark}
There is a partially symmetric analog of the apolarity lemma, which is a special case of \cite[Theorem 1.4]{Maciek}. We recall that the annihilator of the Jacobian ideal $J_F = (y_0[F],\hdots, y_{n}[F]) \subset R$ of $F$ is the ideal
\[
\Ann(J_F)=\{\delta\in S\mid \delta[y_0[F]]=\dots=\delta[y_n[F]]=0\}\subset S
\]
of all differential polynomials annihilating all partial derivatives of $F$.

\begin{lemma}\label{lemma: partially symm apolarity} Let $F\in R_d$ be a form of degree $d$, let $[\ell_1],\dots,[\ell_r]\in \PP^n(R_1)$ be $r$ distinct linear forms and $Z=\{\ell_1^{\vee},\dots,\ell_r^{\vee}\}\subset \PP^n$. Then the partial derivatives $y_0[F],\dots,y_n[F]$ can be simultaneously decomposed as linear combinations of $\ell_1^{d-1},\dots,\ell_r^{d-1}$ if and only if $I(Z)\subset \Ann(J_F)$. 

In particular, $\prk(F)\le r$ if and only if $\Ann(J_F)$ contains the ideal of a set of $r$ points.
\end{lemma}

\begin{remark}\label{remark: the annihilators of F and JF are very similar}
If $F\in R_d$, then $\Ann(F)\subset\Ann(J_F)$ by definition. Moreover, if $\delta\in \Ann(J_F)$ is a differential form of degree $e\in\{1,\dots,d-1\}$, then $0 = \delta[y_i[F]] = y_i[\delta[F]]$ for every $i \in\{ 0,\hdots, n\}$. Hence $\delta[F] = 0$. Thus, $\Ann(F)_{e} = \Ann(J_{F})_e$ for every $e \in\{ 0, \hdots,d-1\}$. However, every differential form of degree $d$ annihilates every form of degree $d-1$, so $\Ann(J_F)_{d}=S_d$, while $\ann(F)_d\subsetneq S_d$. 
\end{remark}

We will take advantage of a few families for which the rank is known. First, the symmetric and partially symmetric rank of a quadratic form $F \in R_2$ both equal the rank of the Hessian matrix of its second derivatives. For binary forms, the symmetric and partially symmetric rank are the same by \cite[Corollary 3.12]{ZHQcomonbinary}, and we can compute them with Sylvester's algorithm. Here we recall a version from 
\cite[Section 3.2]{deParisMaximum}.

\begin{proposition} \label{propos: Sylvester's algorithm}
Let $F\in\CC[x_0,x_1]$ be a binary form. Then $\Ann(F)$ is a complete intersection generated by two elements $\delta_1,\delta_2$. Assume that $\deg(\delta_1)\le \deg(\delta_2)$. Then
\[
\prk(F)=\rk(F)=\begin{cases}
  \deg(\delta_1) & \mbox{ if $\delta_1$ is squarefree}\\
\deg(\delta_2) & \mbox{ otherwise.}
\end{cases}
\]
\end{proposition}

Finally, monomials' symmetric and partially symmetric rank coincide by \cite[Theorem 3.8]{PartComon}. It has been first computed in \cite[Proposition 3.1]{rankofmonomials}.

\begin{proposition}\label{proposition: rank of monomials}
Let $a_0\le \dots\le a_n$ be positive integers. Then
\[
\prk(x_0^{a_0}\cdots x_n^{a_n})=\rk(x_0^{a_0}\dots x_n^{a_n})=\prod_{i=1}^n(1+a_i).
\]
\end{proposition}

\begin{remark}\label{remark: catalecticant bound}
One way to bound the rank of a polynomial is the \emph{catalecticant method}, that states that if $F\in R_d$ and $t\in\{1,\dots,d\}$, then
\[\rk(F)\ge \binom{n+t}{n}-\dim\ann(F)_t.\]
Indeed, if $Z\subset\p^n$ is a set of $\rk(F)$ points such that $I(Z)\subset\ann(F)$, then $\dim I(Z)_t\ge \binom{n+t}{n}-\deg(Z)$, so
\[
\rk(F)= \deg(Z)\ge \binom{n+t}{n}-\dim I(Z)_t\ge \binom{n+t}{n}-\dim \ann(F)_t.
\]
\end{remark}

We end this section with a more geometric lemma we use later in the paper. For the sake of completeness, we include a simple proof. 
\begin{lemma}\label{lemma: can find n-k independent points not spanning L}
Let $n\ge 2$ and let $k\in\{0,\dots,n-1\}$. Let $\cL\subset\PP^n$ be a linear subspace of dimension $k$ and let $Z\subset\PP^n$ be a set of points such that $\langle Z\rangle + \cL =\PP^n$. Then $Z$ contains a subset $S$ of $n-k$ linearly independent points such that $\langle S \rangle \cap \cL =\emptyset$.
\end{lemma}

\begin{proof} Let $S \subset Z$ be a maximal subset of  linearly independent  points such that $\langle S \rangle \cap \cL = \emptyset$. We only need to show that $|S| = n-k$. To this end, it suffices to prove that $Z \subset \langle S \rangle + \cL$. Indeed, then we have $\langle S \rangle + \cL = \langle Z \rangle + \cL = \PP^n$
by hypothesis, and we obtain $\dim \Span{S} = n-k-1$.

Let $q \in Z$. If $q \in \langle S \rangle$, then there is nothing to prove. If $q \notin \langle S \rangle$, then $\{q\}\cup S$ is a linearly independent subset of $Z$. By maximality of $S$, $\langle S \cup \{q\} \rangle \cap \cL = \{q_{\cL}\}$ is a point of $\cL$. We obtain $q \in \langle S \cup \{q_{\cL}\} \rangle \subset \langle S \rangle + \cL$, as wanted.
\end{proof}

\section{A lower bound}
\label{Section: New A lower bound}

In this section, we work on a lower bound for the symmetric rank of a product of linear forms and a quadratic form $F = \ell_0^{a_0}\cdots \ell_{k}^{a_k}\cdot q \in R$, assuming some properties of their annihilators. These assumptions turn out to be instead natural; for instance, they hold when $F$ is a cubic or a quartic, and they are satisfied by some large families in higher degree (see Example \ref{Example: high degrees examples} and Section \ref{Section: The symmetric rank of quartics}). 
We illustrate our bound and we test it with several examples. 

\begin{definition}\label{def: conciseness}
A form $F\in R$ is said to be \emph{concise} if $\Ann(F)_1=0$. In other words, there is no change of coordinates such that $F$ is a form in less than $n+1$ variables. 
\end{definition}

Without loss of generality, we work with concise polynomials in the remainder of the paper. Precisely, let $F = \ell_0^{a_0} \cdots \ell_k^{a_k} \cdot q \in R$ be a concise form of degree $\deg(F) = a_0 + \cdots + a_k + 2$, with $\ell_0,\dots, \ell_k \in R_1$ pairwise  different linear forms, and $q \in R_2$ a quadratic form. We denote $\cL = \PP(\langle \ell_0^{\vee}, \hdots, \ell_{k}^{\vee} \rangle) \subset \PP^n$. We aim to understand the symmetric rank of this type of polynomial. Our result provides a lower bound that depends on $\deg(F)$ and the codimension of the space $\cL$. In many cases the bound is a linear function on $n$. This is not a coincidence but rather a feature of the chosen type of polynomials. We make this observation clear in the next remark.

\begin{remark}\label{rmk: il rango di F è lineare in n}
If we fix the degree $\deg(F)$ and the dimension $\dim \cL$, then asymptotically in $n$, the symmetric rank of $F$ can be bounded by a linear function on $n$. We argue as follows. Set $\dim \cL = j$. Up to $\GL-$equivalence, we can take 
\[F = x_0^{a_0}\cdots x_{j}^{a_j}\ell_{j+1}^{a_{j+1}}\cdots \ell_{k}^{a_k}(g + x_{j+1}^2 + \cdots + x_n^2),\]
with $\ell_{j+1},\hdots,\ell_k, g \in \CC[x_0,\hdots, x_{j}]$. On the one hand, we have
\[\rk(x_0^{a_0}\cdots x_{j}^{a_j}\ell_{j+1}^{a_{j+1}}\cdots \ell_{k}^{a_k}\cdot g) \leq \binom{j + \deg(F)}{j},\]
since if $f \in R_d$, then $\rk(f) \leq \dim R_d$ (see \cite[Section 3.3]{deParisMaximum}). 
On the other hand, the maximum rank of a monomial of degree $d$ in $j+2$ variables is $M=(1 + \frac{d-1}{j+1})^{j+1}$ using Proposition \ref{proposition: rank of monomials}. 
So $\rk(x_0^{a_0}\cdots x_{j}^{a_j}\ell_{j+1}^{a_{j+1}}\cdots \ell_{k}^{a_k}(x_{j+1}^2 + \cdots + x_n^2)) \leq M\binom{\deg(F)+j}{j}(n-j)$ for large $n$. Combining both bounds, we can assure the existence of a constant $K \in \CC$ such that $\rk(F) \leq Kn$ asymptotically in $n$. This behavior is no longer true for arbitrary $\deg(F)$ and $\dim \cL$, as we discuss in Example \ref{Example: non asymtotic in n}. 
\end{remark}

Before continuing, we introduce the following notation. For $H_{1},\hdots,H_{i} \in S_1$, we write the action of $H_{i}\cdots H_{1}$ on $F \in R$ as 
\begin{equation}\label{Equation: action on F} \begin{array}{lllllllllllllll}
H_{i}\cdots H_{1}[F] & = & H_{i}\cdots H_{1}[\ell_0^{a_0} \cdots \ell_k^{a_k}] \cdot q 
+ \sum_{j} \frac{H_{i}\cdots H_{1}}{H_{j}}[\ell_0^{a_0} \cdots \ell_k^{a_k}] \cdot H_{j}[q] \\[0.15cm]
     & & + \sum_{j\neq m} \frac{H_{i}\cdots H_{1}}{H_{j}H_{l}}[\ell_0^{a_0} \cdots \ell_k^{a_k}] \cdot H_{j}H_{m}[q].
\end{array}
\end{equation}
For $H \neq 0$ and $\ell \neq 0$, the vanishing $H[\ell] = 0$ is equivalent to say that the dual point $\ell^{\vee} = (a_0:\cdots :a_n) \in \PP^n$ belongs to the hyperplane $\Var(H) \subset \PP^n$. We sometimes denote the linear form $H$ by $H_c$ to emphasize $c = (c_0:\cdots:c_n) \in \PP^n$.

Our main result is the following (see next Remark \ref{Remark: lower bound extended} for a further generalization). 

\begin{theorem}\label{Proposition: generalized lower bound}
Let $F = \ell_0^{a_0}\cdots \ell_k^{a_k} \cdot q \in R$ be concise  with $\dim \cL = k < n$.  Assume $F$ satisfies the following two conditions. 
\begin{itemize}
    \item[(a)] For every $2\le i \leq \deg(F)-1$ and $H_{i}\cdots H_{1} \in \Ann(F) \subset S$, if $\cL \subset \Var(H_{i})$, then $H_{i-1}\cdots H_{1}[\ell_0^{a_0}\cdots \ell_k^{a_k}] = 0$. 
   \item[(b)] If $\dim \cL > 0$ and $\cL \subset \Var(H)$ for some $H\in\Ann(F)_1$, then $H[q]^{\vee} \notin \cL$. 
\end{itemize}
If $\dim \cL = 0$, then $\rk(F) \geq (\deg(F) - 1)n$.  If $\dim \cL > 0$, then $\rk(F) \geq \deg(F)(n - k)$. 
\end{theorem}

\begin{proof} Assume $Z \subset \PP^n$ is a set of points with $I(Z) \subset \Ann(F) \subset S$. For every $j\in\{0,\dots,k\}$ we set $\alpha_j = a_0 + \cdots + a_j$. In particular $\deg(F) = \alpha_{k}+2$.  

To show the first statement, we will prove that condition (a) implies that $Z$ can be written as a disjoint union $Z = Z_1 \cup \cdots \cup Z_{\alpha_k + 2}$ such that for every $t\in\{1,\dots,\alpha_k+1\}$, the set $Z_t$ consists of $n-k$  linearly independent  points satisfying $\langle Z_t \rangle \cap \cL = \emptyset$. This implies that $|Z| \geq (\deg(F)-1)(n-k)$. As a consequence of Lemma \ref{lemma: apolarity}, $\rk(F) \geq (\deg(F)-1)(n-k)$.

The conciseness of $F$ guarantees that $\langle Z\rangle=\PP^n$, so by Lemma \ref{lemma: can find n-k independent points not spanning L} there exists a subset $Z_1\subset Z$ of $n-k$  linearly independent  points such that $\langle Z_1 \rangle \cap \cL = \emptyset$. If $\langle Z\setminus Z_1\rangle + \cL =\PP^n$, then Lemma \ref{lemma: can find n-k independent points not spanning L} implies the existence of a subset $Z_2\subset Z\setminus Z_1$ of $n-k$  linearly independent  points such that $\langle Z_2 \rangle \cap \cL = \emptyset$. We iterate this process as many times $i \geq 2$ as possible; we find $i-1$ subsets $Z_1, \hdots, Z_{i-1}$, each cnsisting of $n-k$  linearly independent   points such that $\langle Z_t \rangle \cap \cL = \emptyset$ for every $t\in\{1,\dots,i-1\}$; and $Z_i=Z\setminus (Z_1\cup\dots\cup Z_{i-1})$ with $\dim \, \langle Z_i \rangle + \cL < n$. We take a hyperplane $\Var(H_{i}) \subset \PP^n$ containing $\langle Z_i \rangle + \cL$. 
Assume by contradiction that $i\le \alpha_k+1$ and define
\[ j=\min
\{t\in\{0,\dots,k\}\mid i\le \alpha_t+1
\}.\]
Now, we define $i-1$ more hyperplanes as follows. If $j = 0$, we define $i-1$ hyperplanes 
\begin{align*}
\Var(H_1)=\Span{Z_1,\ell_1^\vee,\dots,\ell_k^\vee},\dots,
\Var(H_{i-1})=\Span{Z_{i-1},\ell_1^\vee,\dots,\ell_k^\vee} 
\end{align*}
not containing $\ell_0^{\vee}$. For $j \geq 1$, instead we define $a_0$ hyperplanes 
\begin{align*}
\Var(H_1)=\Span{Z_1,\ell_1^\vee,\dots,\ell_k^\vee},\dots,
\Var(H_{\alpha_0})=\Span{Z_{\alpha_0},\ell_1^\vee,\dots,\ell_k^\vee} 
\end{align*}
not containing $\ell_0^{\vee}$, then we take $a_1$ hyperplanes 
\begin{align*}
\Var(H_{\alpha_0+1})=\Span{Z_{\alpha_0+1},\ell_0^\vee,\ell_2^\vee,
\dots,\ell_k^\vee},
\dots,
\Var(H_{\alpha_1})=\Span{Z_{\alpha_1},\ell_0^\vee,\ell_2^\vee,\dots,\ell_k^\vee} 
\end{align*}
not containing $\ell_1^{\vee}$, and we continue in this way up to defining $a_{j-1}$ hyperplanes 
\begin{align*}
\Var(H_{\alpha_{j-2}+1})=\Span{Z_{\alpha_{j-2}+1},\ell_0^\vee,\dots,\ell_{j-2}^\vee,l_j^\vee,\dots,\ell_k^\vee},
\dots,
\Var(H_{\alpha_{j-1}})&=\Span{Z_{\alpha_{j-1}},\ell_0^\vee,\dots,\ell_{j-2}^\vee,\ell_j^\vee,\dots,\ell_k^\vee}
\end{align*}
not containing $\ell_{j-1}^{\vee}$. Lastly, we choose $i-\alpha_{j-1}-1$ hyperplanes
\begin{align*}
\Var(H_{\alpha_{j-1}+1})=\Span{Z_{\alpha_{j-1}+1},\ell_0^\vee,\dots,\ell_{j-1}^\vee,\ell_{j+1}^\vee,\dots,\ell_k^\vee},
\dots,
\Var(H_{i-1})&=\Span{Z_{i-1},\ell_0^\vee,\dots,\ell_{j-1}^\vee,\ell_{j+1}^\vee,\dots,\ell_k^\vee}
\end{align*}
not containing $\ell_{j}^{\vee}$.

By construction $H_{i}\cdots H_{1} \in I(Z) \subset \Ann(F)$ and $i \leq \deg(F) -  1$, hence $H_{i-1}\cdots H_{1}[\ell_0^{a_0} \cdots \ell_k^{a_k}] = 0$ by hypothesis (a). Moreover, if $s\in\{0,\dots, k\}$ and $t\in\{1,\dots,i-1\}$, then $H_{t}[\ell_s] \neq 0$ if and only if $\alpha_{s-1} +1\leq t \leq \alpha_{s}$. Thus
\begin{equation}\label{Equation: set splitting}
0 = H_{i-1}\cdots H_{1}[\ell_0^{a_0}\cdots \ell_k^{a_k}] = K\cdot(H_{i-1}[\ell_{j}]\cdots H_{1}[\ell_0])(\ell_{j}^{a_j-i+1}\ell_{j+1}^{a_{j+1}} \cdots \ell_k^{a_k})
\end{equation}
for some non-zero constant $K \in \CC$. This implies $0 \neq H_{i-1}[\ell_{j}]\cdots H_{1}[\ell_0] = 0$, a contradiction
.  

Now, we consider the second statement. Assume $\dim \cL > 0$ and write $Z = Z_1 \cup \cdots \cup Z_{\alpha_k+2}$ as above.
Suppose by contradiction that  $|Z| \leq  (\alpha_{k}+2)(n - k)-1$
. Therefore, $|Z_{\alpha_k+2}| \leq n-k-1$ by construction and we can take a hyperplane $\Var(H_{\alpha_k + 2})$ containing $\langle Z_{\alpha_k + 2} \rangle + \cL$. We set $Q= H_{\alpha_k +2}[q]^{\vee}$, which is well defined: indeed, if $H[q] = 0$ and $\cL \subset \Var(H)$, we would have $H[F] = 0$, contradicting the conciseness of $F$. We have $Q \notin \cL$ by (b). We define $\pi = \langle \{Q\} \cup \cL \rangle\subset \PP^n$, a linear subspace of dimension $k+1$. By construction $\langle Z_{\alpha_k+1} \rangle \cap \cL = \emptyset$, hence by Lemma \ref{lemma: can find n-k independent points not spanning L} we can find a subset $Y_{\alpha_k+1} = Z_{\alpha_k+1} \setminus \{P\}$ of $n-k-1$  linearly independent  points such that $\langle Y_{\alpha_k+1} \rangle \cap \pi = \emptyset$. So, we can take a hyperplane $\Var(H_{\alpha_k +1})$ containing $\langle Y_{\alpha_k+2} \rangle + \cL$ and not containing  the point $Q$. Thus $H_{\alpha_k +1}H_{\alpha_k +2}[q] \neq 0$ which is equivalent to $Q \notin \Var(H_{\alpha_k +1})$. 

Now, we set $Y_1 = Z_1 \cup \{P\}$ and we take a hyperplane $\Var(H_{1})$ containing $Y_1$ but not $\cL$. Without loss of generality, we can assume that $\ell_0^{\vee} \notin \Var(H_{1})$, i.e. $H_{1}[\ell_0] \neq 0$. Notice that we are using $\dim \cL > 0$, otherwise $Y_1$ could generate $\PP^n$. In addition, for  $2 \le t \leq \alpha_k$ with $\alpha_{s-1}+1 \leq t \leq \alpha_{s}$, we take hyperplanes $\Var(H_{t})= \langle Z_t \rangle + \langle \ell_0^{\vee},\hdots, \ell_{s-1}^{\vee}, \ell_{s+1}^{\vee}, \hdots, \ell_k^{\vee} \rangle$. By construction, 
$H_{\alpha_{k}+2}\cdots H_{1} \in I(Z) \subset \Ann(F)$ and $\cL \subset \Var(H_{\alpha_k +2})\cap \Var(H_{\alpha_k +1})$. Inserting this in (\ref{Equation: action on F}), we get
$0 = (H_{\alpha_k +2}\cdots H_{1})[F] = 
(H_{\alpha_k}\cdots H_{1})[\ell_0^{a_0} \cdots \ell_k^{a_k}]\cdot (H_{\alpha_{k}+1}H_{\alpha_k+2})[q].$ However, 
\begin{equation}\label{Equation: large splitting}
0 = (H_{\alpha_k}\cdots H_{1})[\ell_0^{a_0}\cdots \ell_k^{a_k}] = 
K\cdot H_{\alpha_k}[\ell_k]\cdots H_{1}[\ell_0] \neq 0,
\end{equation}
for some non-zero constant $K \in \CC$, and we arrive at a contradiction. 
\end{proof}

\begin{remark}\label{Remark: lower bound extended} (i) In the statement of Theorem \ref{Proposition: generalized lower bound}, we assume that $\ell_0,\dots,\ell_k$ are linearly independent, but we can generalize the result. Let $F = \ell_0^{a_0} \cdots \ell_{k}^{a_k} \cdot q$ with  $\dim \cL = j < k$ and $\cL = \langle \ell_0^{\vee} ,\hdots, l_j^{\vee} \rangle$. Set $d = a_0 + \cdots + a_j$. Assume $F$ satisfies (a): for $i \leq d+1$ and $H_i\cdots H_1 \in \Ann(F)$, if $\cL \subset \Var(H_i)$, then $H_{i-1}\cdots H_{1}[\ell_0^{a_0}\cdots \ell_{j}^{a_j}] = 0$. Then  
\[\rk(F) \geq (d+1)(n-j).\]
The arguments are verbatim to those exposed in the first part of the above proof.

\vspace{0.15cm}
\noindent (ii) Let $F = \ell_0^{a_0}\cdots \ell_k^{a_k} \cdot q$ with $\dim \cL = k$. A weaker version of condition (a) in Theorem \ref{Proposition: generalized lower bound} is: there is an integer $2 \leq \Gamma \leq \deg(F) - 1$ such that for all $2 \leq i \leq \Gamma$ and $H_i \cdots H_1 \in \Ann(F)$, if $\cL \subset \Var(H_i)$, then $H_{i-1}\cdots H_i[\ell_0^{a_0}\cdots \ell_{k}^{a_k}] = 0$. If $F$ satisfies it, the same arguments prove that $\rk(F) \geq \Gamma(n-k)$. 
\end{remark}

As we have pointed out in Remark \ref{rmk: il rango di F è lineare in n}, the symmetric rank of $F$ is a linear function on $n$ provided that $\deg(F) \ll n$ and  $\dim \cL \ll n$. Interestingly, Theorem \ref{Proposition: generalized lower bound} gives, as far as we know, a new non-trivial bound for $\rk(F)$. Otherwise, the $\rk(F)$ behavior is barely known. For instance, when $\dim \cL$ is close to $n$, $\rk(F)$ may not be linear on $n$ while Theorem \ref{Proposition: generalized lower bound} provides $\rk(F) \geq K\deg(F)$. When $\deg(F)$ is close to $n^{s}$, for some integer $s$, and $\dim \cL \ll n$, then Theorem \ref{Proposition: generalized lower bound} provides a bound $\rk(F) \geq K n^{s+1}$. Let us illustrate it with some examples.  

\begin{example}\label{Example: non asymtotic in n} (i) For $F = x_0\cdots x_{n}(x_0+\cdots + x_{n})(x_0^2 + \cdots + x_n^2)$ we have $\rk(F) \geq \binom{n+2}{2}$. 
In order to prove that, by Remark \ref{remark: catalecticant bound} it is enough to show that $\ann(F)_2=0$.

Realizing $\Ann(F)_2$ as the kernel of the map $\rho: S_2 \to R_{d-2}$ sending $g \to g[F]$, it is enough to see that $\rho$ has rank $\binom{n+2}{2}$. Consider the usual monomial basis $B_1$ of $S_2$ ordered lexicographically with $y_n \geq \cdots \geq y_0$; and the usual monomial basis $B_2$ of $R_{d-2}$ ordered lexicographically with $x_0 \geq \cdots \geq x_n$. Let $M$ be the matrix associated with $\rho$ in these bases. Set $H = \langle g[F] \mid g \in B_1 \setminus \{y_0y_1, \hdots, y_0y_n\} \rangle$ and $E = \langle y_0y_1[F], \hdots, y_0y_n[F] \rangle$. After permutation, the $\rho_{|H}$ matrix defines an under triangular matrix. We write $y_0y_i[F] = x_1^4x_2 \cdots x_{i-1}x_{i+1} \cdots x_n + f_i$ for $i = 2,\hdots,n$ and $y_0y_1[F] = x_2^4x_3\cdots x_n + f_1$. We have that for $i = 2,\hdots, n$, the monomial $x_1^4x_2 \cdots x_{i-1}x_{i+1}\cdots x_n$ occurs in $g[F]$ if and only if $g = y_0y_i$; and $x_2^4x_3\cdots x_n$ occurs in $g[F]$ if and only if $g = y_0y_1$. Hence, $\rho(R_2) = H \oplus E$ and the result follows. 

\vspace{0.25cm} (ii) Take $F = x_0\cdots x_{k}(x_0+\cdots + x_{k})(x_{k+1}^2 + \cdots + x_n^2)$ with $k < n$. This form satisfies condition (a) for $2 \leq i \leq k+2$ in Theorem \ref{Proposition: generalized lower bound} by Examples \ref{Example: high degrees examples}(i) and (ii), so we have $\rk(F) \geq (k+2)(n-k)$ by Remark \ref{Remark: lower bound extended}(i). Now, using the above example (i), we can compute
\[\Ann(F)_2 = \langle y_{i}y_j \mid k+1 \leq i < j \leq n \rangle + \langle y_{k+1}^2 - y_j^2 \mid j = k+2,\hdots n \rangle,\]
taking into account that $\Ann(F) = \CC[x_{k+1},\hdots,x_n]\Ann(x_0\cdots x_k(x_0+\cdots+x_k))_{|\CC[x_0,\hdots,x_k]} + \CC[x_0,\hdots, x_k]\Ann(x_{k+1}^2 + \cdots + x_n^2)_{|\CC[x_{k+1},\hdots, x_n]}$. 
In particular, $\dim \Ann(F)_2 = \binom{n-k}{2} + n-k-1$. By Remark \ref{remark: catalecticant bound}, we deduce that
\[
\rk(F)\ge \binom{n+2}{2}-\binom{n-k}{2} - n+k+1
\]
For low codimension $n-k$, then this gives $\rk(F) \geq K n^2$, while Theorem \ref{Proposition: generalized lower bound} only provides $\rk(F) \geq  (k+2)(n-k) \sim  n$. 

\vspace{0.25cm} (iii)  Take $F = x_0^{n-3}x_1(x_{2}^2 + \cdots + x_n^2)$ with $n \geq 5$. This form satisfies the hypothesis of Theorem \ref{Proposition: generalized lower bound} by Example \ref{Example: high degrees examples}(i) and (ii), so we have $\rk(F) \geq n(n-1)$. Moreover we can bound $\rk(F) \leq \rk(x_0^{n-3}x_1(x_2^2 + x_3^2)) + \rk(x_0^{n-3}x_1(x_4^3 + x_5^2)) + \cdots \leq 2n(n-2)$, using $\rk(x_0^{n-3}x_1(x_i^2 + x_{i+1}^2)) = \rk(x_0^{n-3}x_1x_ix_{i+1}) = 4(n-2)$ by Proposition \ref{proposition: rank of monomials}. 
\end{example}

The following result gathers two geometric necessary conditions on a set of points $Z\subset \PP^n$ for its ideal $I(Z) \subset S$ to be contained in $\Ann(F)$; and a bound for the partially symmetric rank $\prk(F)$. They follow directly from the proof of Theorem \ref{Proposition: generalized lower bound}. 

\begin{corollary}\label{Corollary: geometric conse generalized lower bound} Let $F = \ell_0^{a_0}\cdots \ell_{k}^{a_k} \cdot q$ be a concise form of degree $\deg(F)$ with $\dim \cL = k$. Assume $F$ satisfies condition (a) in Theorem \ref{Proposition: generalized lower bound}. 
\begin{itemize}
    \item [(i)] If $I(Z) \subset \Ann(F)$, then $Z = Z_1 \cup \cdots \cup Z_{\deg(F)}$ with $Z_i$ a set of $n- \dim \cL$  linearly independent  points such that $\langle Z_i \rangle \cap \cL = \emptyset$ for every $i\in\{ 1,\hdots, \deg(F) - 1\}$.
    Consequently, $\rk(F) \geq (\deg(F)-1)(n - \dim \cL)$.
    \item[(ii)] Suppose $\dim \cL > 0$ and let $Z \subset \PP^n$ be a set of points admitting a splitting $Z = Z_1 \cup \cdots \cup Z_{\deg(F)}$ as in (i). If there is a hyperplane $\Var(H)$ containing $\langle Z_{\deg(F)} \rangle + \cL$ and such that $Q=H[q]^{\vee} \notin \cL$, then $I(Z) \not\subset \Ann(F)$. 
     \item [(iii)] $\prk(F) \geq (\deg(F)-1)(n - \dim \cL)$.
\end{itemize}
\end{corollary}
\begin{proof} (iii) Using Remark \ref{remark: the annihilators of F and JF are very similar}, we have that $\Ann(J_{F})_{i} = \Ann(F)_i$ for every $0 \leq i \leq \deg(F)-1$. Therefore, the same arguments proving (i) in the proof of Theorem \ref{Proposition: generalized lower bound} also shows that if $I(Z) \subset \Ann(J_{F})$, then $Z = Z_1 \cup \dots \cup Z_{\deg(F)}$ with $Z_i$ a set of $n- \dim \cL$  linearly independent  points such that $\langle Z_i \rangle \cap \cL = \emptyset$ for every $i\in\{ 1,\hdots, \deg(F) - 1\}$. Consequently, $\prk(F) \geq (\deg(F) -1)(n- \dim \cL)$. 
\end{proof}

We apply this lower bound to determine the symmetric and \psr of the product of a high multiplicity hyperplane and a quadric
. If the product of an irreducible quadric and a multiplicity $d-2$ hyperplane is concise, then it is $\GL-$equivalent to one of the forms (i), (ii), or (iii) below. In the special case $d = 3$, the symmetric rank is computed in \cite[Proposition 7.2]{Landsberg-Teitler} and \cite[Theorem 4.5]{reduciblecubics}.

\begin{proposition}\label{Proposition: rank of high multiplicity hyperplane}
Let $n \geq 2$ and $d \geq 3$. Then
\begin{itemize}
\item[(i)] $\rk(x_0^{d-2}(x_0^2+\dots+x_n^2)) = \prk(x_0^{d-2}(x_0^2+\dots+x_n^2)) = (d-1)n$.
\item[(ii)] $\rk(x_0^{d-2}(x_1^2+\dots+x_n^2)) = \prk(x_0^{d-2}(x_1^2+\dots+x_n^2)) = (d-1)n$.
\item[(iii)] $\rk(x_0^{d-2}(x_0x_1+x_2^2+\dots+x_n^2)) = \prk(x_0^{d-2}(x_0x_1+x_2^2+\dots+x_n^2)) = (d-1)n + 1$.
\end{itemize}
\end{proposition}
\begin{proof} For convenience, we denote the forms by $F_1, F_2$ and $F_3$ and $Q_0 = (1: 0: \cdots :0) \in \PP^n$. For the lower bound, we want to apply Corollary \ref{Corollary: geometric conse generalized lower bound} with $\cL = \{Q_0\}$. Fix $i \leq d - 1$ and a product of linear forms $H_{i}\cdots H_{1} \in \Ann(F_k) \subset S$ with $Q_0 \in \Var(H_{i})$. We write $H_{j} = c_0^jy_0 + \cdots + c_n^jy_n$. 

For $F \in \{F_1,F_2,F_3\}$, from (\ref{Equation: action on F}) we have
\[\begin{array}{lllllllllllllll}
H_{i}\cdots H_{1}[F]  =  H_{i-1}\cdots H_{1}[x_0^{d-2}] \cdot H_{i}[q] + \sum_{j < i} \frac{H_{i}\cdots H_{1}}{H_{j}H_{i}}[x_0^{d-2}] \cdot H_{j}H_{i}[q], 
\end{array}\]
formally a polynomial of degree $d-i$. If $j \in \{1,\hdots,n\}$, then the coefficient of $x_0^{d-i-1}x_j$ is
\[\frac{(d-2)!}{(d-i-1)!} c^i_{j}c_0^{i-1}\cdots c_0^{1}.\]
Since $H_{i}\cdots H_{1}[F] = 0$ and $c^i_j \neq 0$ for some $j \in \{1,\hdots,n\}$, we have $c_0^{i-1}\cdots c_0^1 = 0$. Thus, $H_{i-1}\cdots H_{1}[x_0^{d-2}] = 0$. 
Remark \ref{remark: sym rank bigger than par sym rank} and Corollary \ref{Corollary: geometric conse generalized lower bound}(iii) assure that 
\[
\prk(F) \geq \rk(F) \geq (d-1)n.\]

Now, we improve the lower bound for $F_3$. If $H_{d-1}\cdots H_{1}\in \Ann(F_3)$, then $H_{d-1}\cdots H_{1}[F_3]$ is formally a linear form and the coefficient of $x_1$ is $(d-1)!c_0^{d-1}\cdots c_0^{1}$. We can conclude that if $H_{d-1}\cdots H_{1} \in \Ann(F_3)$, then $(d-1)!c_0^{d-1}\cdots c_0^1 = 0$, so at least one hyperplane $\Var(H_{j})$ contains the point $Q_0$. Assume by contradiction that there exists a set $Z \subset \PP^n$ of $(d-1)n$ points such that $I(Z) \subset \Ann(F_3)$. Then $Z = Z_1 \cup \cdots \cup Z_{d-1}$ where $Z_i$ is a set of $n$  linearly independent  points verifying $Q_0 \notin \langle Z_i \rangle$ for any $i \in\{ 1,\hdots, d-1\}$ by Corollary \ref{Corollary: geometric conse generalized lower bound}(i). Taking hyperplanes $\Var(H_{1}) = \langle Z_i \rangle, \dots, \Var(H_{d-1}) = \langle Z_{d-1} \rangle$,  we have $H_{1}\cdots H_{d-1} \in I(Z) \subset \Ann(F_3)$. By construction, $Q_0 \notin \Var(H_{i})$ for any $i\in\{ 1,\hdots, d-1\}$ and we have arrived at a contradiction. Thus $\rk(F_3) \geq (d-1)n + 1$. Notice that the argument above only involves forms of degree at most $d-1$. By Remark \ref{remark: the annihilators of F and JF are very similar} and Lemma \ref{lemma: partially symm apolarity}, the same proof shows that $\prk(F_3) \geq (d-1)n+1$. 

Next we show the upper bounds. The result is true for $d = 3$ by \cite[Proposition 4.3]{reduciblecubics}. Assume $d \geq 4$. Either Proposition \ref{propos: Sylvester's algorithm} or \ref{proposition: rank of monomials} show that $\rk(x_0^{d-2}x_i^2) = d-1$ for every $i \neq 0$, so
\begin{equation}\label{Equation: upper bound F1}
\rk(F_1) = \rk(x_0^{d-2}(x_1^2 + \cdots + x_n^2)) \leq \sum_{i=1}^n \rk(x_0^{d-2}x_i^2) = (d-1)n.
\end{equation}
Now we write $\rk(F_2) = \rk(x_0^{d-2}(x_0^2 + \dots + x_n^2)) 
\leq \rk(x_0^{d-2}(x_0^2 + x_1^2)) +\rk(x_0^{d-2}(x_2^2+\dots+ x_n^2))$. We have $\rk(x_0^{d-2}(x_0^2 + x_1^2)) = d-1$ by Proposition \ref{propos: Sylvester's algorithm}, and using (\ref{Equation: upper bound F1}) we get $\rk(F_2) \leq  d-1 + (d-1)(n-1) = (d-1)n.$ Finally, $\rk(F_3) \leq  \rk(x_0^{d-1}x_1) + \rk(x_0^{d-2}(x_2^2 + \cdots + x_n^2))$ and we obtain analogously $\rk(F_3) \leq d+(d-1)(n-1) = (d-1)n+1$. \end{proof}

\begin{remark}\label{remark: examples of supergeneric and maximal rank}
Observe that Proposition \ref{Proposition: rank of high multiplicity hyperplane} gives an alternative, more geometric proof of \cite[Propositions 4.4 and 4.6]{reduciblecubics}. Moreover, it gives us a way to exhibit examples of generic, supergeneric and maximal symmetric rank. Indeed, 
we observe that 
\begin{itemize}
\item for $d=3$, the symmetric ranks of $F_1$ and $F_2$ equal the general rank for $n\in\{2, 4,5,6\}$ and it exceeds it for $n=3$. If we look at \cite{deParisMaximum}, we notice that $\rk(F_3)$ is actually the maximal symmetric rank of a cubic for $n\in\{2,3\}$ and it exceeds the generic rank for $n\in\{2,\dots,6\}$.
\item for $d=4$, the symmetric ranks of $F_1$ and $F_2$ equal the general rank for $n=2$. Moreover $\rk(F_3)$ is the maximal symmetric rank of a ternary quartic, and it equals the generic rank for $n=3$.
\item if $d=5$ and $n=2$, then $F_1$, $F_2$ and $F_3$ have supergeneric rank.
\item if $d=6$ and $n=2$, then $F_1$ and $F_2$ have generic rank, while $F_3$ has supergeneric rank.
\end{itemize}
\end{remark}

We finish this section with more applications of Theorem \ref{Proposition: generalized lower bound}.

\begin{example}\label{Example: high degrees examples} (i) Take $F = x_0\cdots x_k(x_{k+1}^2 + \cdots + x_n^2)$ with $0 < k < n$. Let us check that $F$ satisfies the hypotheses of Theorem \ref{Proposition: generalized lower bound}. Condition (b) is trivially satisfied, so consider (a). Here $\deg(F) = k+3$ and $\cL = \langle x_0^{\vee}, \hdots, x_{k}^{\vee} \rangle =\PP^k$. If $i \leq k+2$ and $H_{i}\cdots H_{1} \in \Ann(F)$ with $\cL \subset \Var(H_{i})$, then 
\begin{equation*}
H_{i}\cdots H_{1}[F] = H_{i-1}\cdots H_{1}[x_0\cdots x_k] \cdot H_{i}[q] + \sum_{j\neq i} \frac{H_{i-1}\cdots H_{1}}{H_{j}}[x_0\cdots x_k] \cdot H_{j}H_{i}[q] = 0.
\end{equation*}
Writing this expression in $\CC[x_0,\hdots, x_k][x_{k+1},\hdots, x_{n}]$, we obtain formally a polynomial of degree one in the variables $x_{k+1}, \hdots, x_n$ with coefficients $c_{j}^{i}H_{i-1}\cdots H_{1}[x_0\cdots x_k]$, respectively. Hence, $c_{j}^{i}H_{i-1}\cdots H_{1}[x_0\cdots x_k] = 0$ for $j = k+1,\hdots,n$. Since $c_{j}^i \neq 0$ for some $j \in \{k+1,\hdots,n\}$, we get $H_{i-1}\cdots H_{1}[x_0\cdots x_k] = 0$. Thus, $F$ satisfies the hypotheses of Theorem \ref{Proposition: generalized lower bound}, and we can conclude that $\rk(F) \geq (k+3)(n-k)$. 

Let us compare this lower bound to the catalecticant bound from Remark \ref{remark: catalecticant bound}. 
It is easy to show that 
\[\Ann(F) = (y_0^2, \hdots, y_k^2) + (y_{k+1}^2 - y_j^2 \mid j = k+2, \hdots, n) + (y_jy_l \mid k+1 \leq j < l \leq n)\] 
hence $\dim\Ann(F)_2 = n + \binom{n-k}{2}$. 
The catalecticant bound is $\rk(F) \geq \binom{n+2}{2}-n - \binom{n-k}{2}$.
The difference between the two bounds is $2n - \frac{k(k+5)}{2}-1$. Hence, for a fixed degree $k+3$, the bound $(k+3)(n-k)$ is asymptotically better than the catalectivant bound.

\vspace{0.15cm}
\noindent (ii) More generally, if $g\in\CC[x_0,\dots,x_k]_2$ is a quadratic form such that $F = (x_0^{a_0} \cdots x_k^{a_k})(g(x_0,\hdots,x_k) + x_{k+1}^2 + \cdots + x_n^2)$ is concise, then the same arguments as the above example (i) give $\rk(F) \geq (a_0 + \cdots + a_k + 2)(n-k)$. Notice that $\rk(F) \leq \binom{\deg(F) + k}{k} + (\prod_{i=0}^{k} 3(a_i+1))(n-k)$. Hence, there is a non-zero constant $K \in \CC$ such that $\rk(F) \leq Kn$ asymptotically in $n$. The precise value of $\rk(F)$ is unknown as far as we know.

\vspace{0.15cm}
\noindent (iii) Take the quartic form $F = x_0x_1(x_2^2 + \cdots + x_n^2)$. We have $\rk(F) \geq 4(n-1)$ from (ii); and we will prove in the next section that $\rk(F) = 4(n-1)$. Thus, the lower bound in Theorem \ref{Proposition: generalized lower bound} is optimal for this quartic. 
\end{example}

\section{Symmetric rank of some reducible quartics}
\label{Section: The symmetric rank of quartics}
Determining the symmetric rank of an arbitrary product of linear forms and a quadratic form is a challenging problem. For reducible cubics, the first step of the argument contained in  \cite{reduciblecubics} is to determine normal forms, up to $\GL$ equivalence. 
In this section we consider reducbible quartics of the form $\ell_0\ell_1 \cdot q$, and we start by producing a list of simple $\GL-$equivalent forms, according to the relative position 
of the two hyperplanes and the quadric. 
Since the rank of ternary quartics has been studied in detail in \cite[Chapter 3]{kleppe}, we focus on $n\ge 3$.


\begin{proposition}\label{Prop: normal forms of reducible quartics}
 Let $n\ge 3$ and $F = \ell_0 \ell_1 \cdot q \in R_4$ be the product of two linear forms $\ell_0,\ell_1 \in R_1$ and an irreducible quadratic form $q \in R_2$. If $F$ is concise, then $F$ is  $\GL-$equivalent to one of the following:
\begin{enumerate}
\item $F_1=x_0^2(x_0^2+\cdots+x_n^2)$, a full rank quadric with a double secant hyperplane.
\item $F_2=x_0^2(x_1^2+\cdots+x_n^2)$, a rank $n$ quadric with a double secant hyperplane.
\item $F_3=x_0^2(x_0x_1+x_2^2+\cdots+x_n^2)$, a full rank quadric with a double tangent hyperplane.
\item $F_4=x_0x_1(x_2^2+\cdots+x_n^2)$, a rank $n-1$ quadric with two secant hyperplanes.
\item $F_5=x_0x_1(x_1x_2+x_3^2+\cdots+x_n^2)$, a rank $n$ quadric, a hyperplane tangent to the quadric and a secant hyperplane not containing the vertex.
\item $F_6=x_0x_1(x_1^2+\cdots+x_n^2)$, a rank $n$ quadric with two secant hyperplanes, one containing the vertex.
\item $F_{7} = x_0(x_0+x_1)(x_1^2 + \cdots + x_n^2)$, a rank $n$ quadric and two secant hyperplanes not containing the vertex.
\item $F_{8,\lambda} = x_0(\lambda x_0+x_1)(x_0^2+\cdots+x_n^2)$, $\lambda\in\CC$, a full rank quadric with two secant hyperplanes. 
\item $F_9=x_0x_1(x_0^2 + x_1x_2 + x_3^2 + \cdots + x_n^2)$, a full rank quadric with a secant hyperplane and a tangent hyperplane.
\item $F_{10,\lambda}=x_0(-\lambda^2 x_0+x_1+\lambda x_2)(x_0x_1+x_2^2+\cdots+x_n^2)$, $\lambda \in\CC$, a full rank quadric with two tangent hyperplanes.
\end{enumerate}   
\end{proposition}
\begin{proof} Write $F = \ell_0\ell_1 \cdot q$, with $\ell_0=a_0x_0+\cdots+a_nx_n$ and $\ell_1=b_0x_0+\cdots+b_nx_n$. We denote by $V \subset \PP^n$ the vertex of the quadric $\Var(q)$. First, we observe $n-1 \leq \rk(q) \leq n+1$. Indeed, if $\rk(q)\le n-2$, then 
$F$ is $\GL-$equivalent to $\ell_0\ell_1(x_3^2 + \cdots + x_n^2)$. In this case,  
$\dim \langle \ell_0,\ell_1,x_3,\hdots,x_n \rangle \le n$, 
contradicting that $F$ is concise. We distinguish three cases depending on $\rk(q)$.
 
\vspace{0.15cm}
\noindent \underline{Case 1}: $\rk(q)=n-1$. Up to $\GL$-action, we can assume that $F = \ell_0\ell_1(x_2^2 + \cdots + x_n^2)$. Since $F$ is concise, we have $\dim \langle \ell_0,\ell_1, x_2,\hdots,x_n \rangle = n+1$. Hence, $F$ is $\GL$-equivalent to $F_4 = x_0x_1(x_2^2 + \cdots + x_n^2)$. 

\vspace{0.15cm}
\noindent \underline{Case 2}: $\rk(q)=n$. The vertex $V$ of $\Var(q)$ is a point that belongs to every tangent space to $\Var(q)$. Moreover $V \notin \Var(\ell_0) \cap \Var(\ell_1)$ because $F$ is concise. We distinguish two subcases depending on $\Var(\ell_0)$ and $\Var(\ell_1)$. 

\vspace{0.15cm}
\noindent \underline{Subcase 2.1}: $\Var(\ell_1)$ is tangent to $\Var(q)$, so $\Var(\ell_0)$ is secant to $\Var(q)$ and $V \notin \Var(\ell_0)$. Notice that $\Var(\ell_0) \neq \Var(\ell_1)$. Up to $\GL$-action, we may assume $F = x_0\ell_1(x_1x_2 + x_3^2 + \cdots+ x_n^2)$. The cubic $\Var(x_0,\ell_1q)$ defines a reducible cubic in $\Var(x_0)=\PP^{n-1}$ consisting of a smooth quadric and a tangent hyperplane. By \cite[Lemma 4.1]{reduciblecubics}, there is an automorphism $\phi$ of $\PP^{n-1}$ sending this cubic to $\Var(x_1(x_1x_2+x_3^2+\cdots+x_n^2))$. Now, we apply the change of coordinates in $\PP^n$ induced by $\phi$ and fixing $x_0$. We obtain $q=x_1x_2+x_3^2+\cdots+x_n^2$ and $\Var(x_0)\cap \Var(\ell_1)=\Var(x_0,x_1)$. Since $\Var(x_0)\neq \Var(\ell_1)$, we can write $\ell_1=b_0 x_0+x_1$. Imposing that $\Var(\ell_1)$ is tangent to $\Var(q)$, we get $\ell_1 = x_1$. In this case, $F$ is $\GL$-equivalent to $F_5 = x_0x_1(x_1x_2 + x_3^2 + \cdots + x_n^2)$. 

 \vspace{0.15cm}
\noindent \underline{Subcase 2.2}: $\Var(\ell_0)$ and $\Var(\ell_1)$ are both secant to $\Var(q)$. Up to $\GL$-action, we may assume that $\ell_0 = x_0$ and $q = x_1^2 + \cdots + x_n^2$. 
If $\Var(\ell_1) = \Var(x_0)$, then $F$ is $\GL$-equivalent to $F_2 = x_0^2(x_1^2 + \cdots + x_n^2)$. Assume $\Var(\ell_1) \neq \Var(x_0)$. 
Analogously, as in Subcase 2.1, we obtain $q = x_1^2+\cdots+x_n^2$ and 
$\Var(x_0)\cap \Var(\ell_1)=\Var(x_0,x_1)$. We write $\ell_1 = b_0x_0 + x_1$. If $V\in \Var(\ell_1)$, then $\ell_1=x_1$ and we obtain that $F$ is $\GL$-equivalent to $F_6 = x_0x_1(x_1^2 + \cdots + x_n^2)$. If $V \notin \Var(\ell_1)$, then $b_0\neq 0$ and a straightforward change of coordinates shows that $F$ is $\GL$-equivalent to $F_7 = x_0(x_0+x_1)(x_1^2 + \cdots + x_n^2)$. 

\vspace{0.15cm}
\noindent \underline{Case 3}: $\rk(q) = n+1$. We distinguish three subcases depending on the hyperplanes $\Var(\ell_0)$ and $\Var(\ell_1)$. 

\vspace{0.15cm}
\noindent \underline{Subcase 3.1}: $\Var(\ell_0)$ and $\Var(\ell_1)$ are both secant to $\Var(q)$. Up to $\GL$-action, we may assume that $F = x_0\ell_1(x_0^2 + \cdots + x_n^2)$.  If $\Var(\ell_1) = \Var(x_0)$, then $F$ is $\GL$-equivalent to $F_1$. Assume $\Var(\ell_1) \neq \Var(x_0)$. The cubic $\Var(x_0, \ell_1q)$ defines a reducible cubic in $\Var(x_0) = \PP^{n-1}$ consisting of a smooth quadric and a secant hyperplane. By \cite[Lemma 4.1]{reduciblecubics}, there is an automorphism $\phi$ of $\PP^{n-1}$ sending the cubic to $\Var(x_1(x_1^2 + \cdots +x_n^2))$. Considering the change of coordinates induced by $\phi$ and fixing $x_0$, we obtain $q = x_0^2 + \cdots + x_n^2$ and $\Var(x_0) \cap \Var(\ell_1) = \Var(x_0,x_1)$. Hence, $\ell_1 = b_0x_0 + x_1$ and $F$ is $\GL$-invariant to the one of the forms $F_{8,\lambda} = x_0(\lambda x_0 + x_1)(x_0^2 + \cdots + x_n^2)$.

\vspace{0.15cm}
\noindent \underline{Subcase 3.2}: $\Var(\ell_0)$ is secant to $\Var(q)$ and $\Var(\ell_1)$ is tangent to $\Var(q)$. Notice that $\Var(\ell_0) \neq \Var(\ell_1)$. Up to $\GL$-action, we may assume that $F = x_0\ell_1(x_0^2 + x_1x_2 + x_3^2 + \cdots + x_n^2)$. The cubic $\Var(x_0, \ell_1q)$ defines a reducible cubic in $\Var(x_0) = \PP^{n-1}$ consisting of a rank $n$ quadratic form and a tangent hyperplane. By \cite[Lemma 4.1]{reduciblecubics}, there is an automorphism $\phi$ of $\PP^{n-1}$ sending the cubic to $x_1(x_1x_2 + x_3^2 + \cdots + x_n^2)$. Considering the change of coordinates induced by $\phi$ and fixing $x_0$, we get $q = x_0^2 + x_1x_2 + x_3^2 + \cdots + x_n^2$ and $\Var(x_0) \cap \Var(\ell_1) = \Var(x_0,x_1)$. We write $\ell_1 = b_0x_0 + x_1$. Imposing that $\Var(\ell_1)$ is tangent to $\Var(q)$, i.e. 
\[\rk(q_{|\ell_1}) = \rk(x_0^2 + b_0x_0x_2 + x_3^2 + \cdots + x_n^2) \leq n-1,\]
we have $b_0 = 0$. Hence, $F$ is $\GL$-invariant to the form $F_9$.

\vspace{0.15cm}
\noindent \underline{Subcase 3.3}: $\Var(\ell_0)$ and $\Var(\ell_1)$ are both tangent to $\Var(q)$. Up to $\GL$-action, we may assume that $F = x_0\ell_1(x_0x_1 + x_2^2 + \cdots + x_n^2)$. If $\Var(\ell_1) = \Var(x_0)$, then $F$ is $\GL$-equivalent to $F_3$. Assume $\Var(\ell_1) \neq \Var(x_0)$. Now, 
if $\ell_1 \in \CC[x_0,x_1]$, then imposing that $\Var(\ell_1)$ is tangent to $\Var(q)$ we obtain $\ell_1 = x_1$, i.e. $\lambda = 0$. 

Else, $\ell_1 \notin \CC[x_0,x_1]$. The cubic $\Var(x_0,x_1,\ell_1q)$ defines a reducible cubic in $\PP^{n-2}$ consisting of a smooth quadric and a secant hyperplane. Indeed, $\Var(q)$ is smooth, $\Var(\ell_1) \neq \Var(x_0)$, and so $\Var(\ell_1) \cap \Var(q) \notin \Var(x_0,x_1) \subset \Var(x_0)$. By \cite[Lemma 4.1]{reduciblecubics}, there is an automorphism $\phi$ of $\PP^{n-2}$ sending the cubic to $\Var(x_2(x_2^2 + \cdots + x_n^2))$.  Applying the induced change of coordinates in $\PP^n$ and fixing $x_0$ and $x_1$, we obtain  $q=x_0x_1+x_2^2+\cdots+x_n^2$ and $\Var(x_0,x_1) \cap \Var(\ell_1) = \Var(x_0,x_1,x_2)$. We write $\ell_1=b_0x_0+ b_1x_1+ b_2x_2$. Imposing that $\Var(\ell_1)$ is tangent to $\Var(q)$, we have $b_1 \neq 0$ and 
\[\rk(q_{|\ell_1})=\rk(-b_0 x_0^2-b_2 x_0x_2+x_2^2+\cdots+x_n^2)\le n-1.\]
Equivalently, $b_1 \neq 0$ and \[\det\begin{pmatrix}
-b_0 & -b_2 &0  &\dots &0\\
-b_2 & 1 &  0 &\dots&0\\
0 & 0 &1 & \ddots&\vdots\\
\vdots & \vdots & \ddots &\ddots &0\\
 0 &0&\dots&0&1
\end{pmatrix}=0\Leftrightarrow b_0+b_2^2=0.\]
Now, $F$ is $\GL$-equivalent to $F_{10,\lambda} = 
x_0(-\lambda^2 x_0 + x_1 + \lambda x_2)(x_0x_1 + x_2^2 + \cdots + x_n^2)$. 
\end{proof}

\begin{remark} (i) Proposition \ref{Prop: normal forms of reducible quartics} is no longer true for $n=2$. Indeed, consider $q=x_0^2+x_1^2-x_2^2$ and let $\ell=x_1-x_2$, a tangent line. Both $\Var(x_0)$ and $\Var(x_1)$ are secant to $\Var(q)$. However, $x_0\ell q$ is not $\GL-$equivalent to $x_1\ell q$ because $\Var(x_0)\cap \Var(\ell)=(0:1:1)\in \Var(q)$, while $\Var(x_1)\cap \Var(\ell)=(1:0:0)\notin \Var(q)$. Similarly, both $x_0x_1q$ and $x_0(x_1-x_2-x_0)q$ give a smooth conic and two secant lines, but $\Var(x_0,x_1,q)= \emptyset$ while $\Var(x_0,x_1-x_0-x_2,q)=(0:1:1)$. When $n\ge 3$, these differences disappear, because if $\dim \Var(\ell_0,\ell_1,q)\ge n-2$ then both $\Var(\ell_0)$ and $\Var(\ell_1)$ are tangent to $\Var(q)$.

\vspace{0.15cm}
\noindent (ii) The above proposition draws an inductive way to find normal forms in higher degrees. For instance, to determine normal forms for quintics $F = \ell_0\ell_1\ell_2 \cdot q$, one can sort out $q$ depending on $\rk(q) \geq n-3$ and its tangential relative position with respect to $\ell_0 = x_0$ and apply the above proposition to $\Var(x_0, \ell_1\ell_2 \cdot q)$. 
\end{remark}

To the reader interested in computing their own examples, a simple Macaulay2's routine (\cite{Macaulay2}) is the following. Given a form $F = \ell_0\ell_1 \cdot q \in R$, we want to determine $\rk(q)$ and the tangential relative positions between $\Var(q)$ and the hyperplanes $\Var(\ell_0)$ and $\Var(\ell_1)$, respectively. For the symmetric rank of a quadratic form, we can compute its Hessian matrix using \texttt{diff} and determine its rank. On the other hand, a hyperplane $\Var(\ell)$ and an irreducible quadric $\Var(Q)$ are tangent if and only if $\rk(Q_{\ell}) \leq \rk(Q)-2$, where $\Var(Q_{\ell})$ is the restriction of the quadric $\Var(Q)$ to the hyperplane $\Var(\ell)$. We can compute $Q_{\ell}$ using the command \texttt{sub}. For example, take $F = 
x_0(x_0+x_1)(x_0^2 + x_0x_2 + x_1x_2 + x_3^2 + x_4^2) \in \CC[x_0,x_1,x_2,x_3,x_4]$. 

\begin{verbatim}
R = QQ[x_0..x_4]
q = x_0^2 + x_0*x_2 + x_1*x_2 + x_3^2 + x_4^2
v = matrix{{x_0,x_1,x_2,x_3,x_4}}
Hessianq = diff(v ** transpose v, q)
rank Hessianq
q0 = sub(q,x_0=>0)
q1 = sub(q,x_1=>-x_0)
Hessianq0 = diff(v ** transpose v, q0)
Hessianq1 = diff(v ** transpose v, q1)
rank Hessianq0
rank Hessianq1
\end{verbatim}
We obtain that $q$ has rank 5, while $q_{\ell_0}$ has rank 4 and $q_{\ell_1}$ has rank 3. Hence $F$ is $\GL-$equivalent to $F_9$, a full rank quadric $\Var(q)$ with a secant hyperplane $\Var(x_0)$ and a tangent hyperplane $\Var(x_0+x_1)$
. 

\subsection{A lower bound for the symmetric rank}
\label{Section: Lower bound}

In this subsection, we work on a lower bound for the symmetric rank of the quartic forms $F_i$. The main result proves that $\rk(F_i) \geq 4(n-1)$ for $i = 4,\hdots,10$ (Proposition \ref{Theorem: lower rank of quartics}) and is based on the results of Section \ref{Section: New A lower bound}. 

\vspace{0.25cm}
We start pointing out that $\rk(F_1) = \rk(F_2) = 3(n-1)$ and $\rk(F_3) = 3(n-1)+1$ by Proposition \ref{Proposition: rank of high multiplicity hyperplane}, and it coincides with their partially symmetric rank. Hitherto, we focus on the remainder of the quartics. We recall that all forms $F_4,\dots,F_{10}$ can be written as $F_i = x_0\ell_i q_i \in R_4$, where $\ell_i = h_0^ix_0 + x_1 + h_2^ix_2 \in \CC[x_0,x_1,x_2]_1$ and $q_i = g_i(x_0,x_1,x_2) + x_3^2 + \cdots + x_n^2 \in R_2$. We denote by $\cL_i = \Var(h_2^iy_1-y_2,y_3,\hdots,y_n) \subset \PP^n$ the line spanned by the dual points $x_0^{\vee} = (1:0:\cdots:0)$ and $\ell_i^{\vee} = (h_0^i:1:h_2^i:0:\cdots:0)$.

The main result is the following. 
\begin{proposition}\label{Theorem: lower rank of quartics} Fix $i \in \{4,\hdots,10\}$ and a reducible quartic form $F_i = x_0\ell_i q \in R_4$. Then $\rk(F_i) \geq 4(n-1)$ and $\prk(F_i) \geq 3(n-1)$.
\end{proposition}
\begin{proof}
To begin with, notice that quartics $F_{i}$ with $i \notin\{ 5,9,10\}$ are already covered by Example \ref{Example: high degrees examples}(ii). 
We want to apply Theorem \ref{Proposition: generalized lower bound} and Corollary \ref{Corollary: geometric conse generalized lower bound} for the remaining quartics.  

First, we see that $F_{10,\lambda} = x_0(-\lambda^2x_0 + x_1 + \lambda x_2)(x_0x_1 + x_2^2 + \cdots + x_n^2)$ satisfies the hypotheses of Theorem \ref{Proposition: generalized lower bound}; and that $F_5 = x_0x_1(x_1x_2+x_3^2 + \cdots +x_n^2)$ and $F_9
= x_0x_1(x_0^2 + x_1x_2 + x_3^2 + \cdots +x_n^2)$ verify condition (a) in Theorem \ref{Proposition: generalized lower bound}. We obtain $\rk(F_{10}) \geq 4(n-1)$ and $\prk(F_i) \geq 3(n-1)$ by Corollary \ref{Corollary: geometric conse generalized lower bound}(iii). 

Let us begin with the hypothesis (a), that is, for $j\in\{2,3\}$ and $H_{j}\cdots H_{1} \in \Ann(F_i)$, if $\cL_i \subset \Var(H_{1})$, then $H_{j}H_{2}[x_0\ell_i] = 0$. We take $F = F_5$, analogous arguments work for $F_9$ and $F_{10,\lambda}$, $\lambda \in \CC$, as well. We have $\cL = 
\Var(y_2,\hdots,y_n)$ and assume $\cL \subset \Var(H_{1})$. Writing $H_{j}\cdots H_{1}[F]$ as a polynomial in 
$\CC[x_0,x_1,x_2][x_3,\hdots,x_n]$, the hypothesis $H_{j}\cdots H_{1}[F] = 0$ gives $H_{1} = y_2$. Moreover, we have $H_{1}[F] = x_0x_1^2$. If $H_{2}H_{1} \in \Ann(F_5)$, we obtain 
\[H_{2}[x_0x_1^2] = H_{2}[x_0]x_1^2 + 2H_{2}[x_1]x_0x_1 = 0,\]
implying $\cL \subset \Var(H_{2})$, so $H_{2}[x_0x_1] = 0$. If $H_{3}H_{2}H_{1} \in \Ann(F_5)$, then 
\[H_{3}H_{2}[x_0x_1^2] = 2(H_{3}[x_1]H_{2}[x_0] + H_{3}[x_0]H_{2}[x_1])x_1 + H_{3}[x_1]H_{2}[x_1]x_0 = 0.\]
Hence, $H_{3}[x_1]H_{2}[x_1] = 0$ and $H_{3}[x_1]H_{2}[x_0] + H_{3}[x_0]H_{2}[x_1] = 0$. Now, $H_{3}H_{2}[x_0x_1] = H_{3}[x_0]H_{2}[x_1] + H_{3}[x_1]H_{2}[x_0] = 0$, as wanted. On the other hand, $F_{10,\lambda}$ also verifies condition (b) in Theorem  \ref{Proposition: generalized lower bound}, for any $\lambda \in \CC$. We have $\cL = \Var(-\lambda x_1 + x_2, x_3,\hdots,x_n)$. If $\cL \subset \Var(H_c)$, then $H_c[q] = c_1x_0 + 2c_2x_2 + \cdots + 2c_nx_n$ and $H_c[q]^{\vee} = (c_1: 0: 2c_2: \cdots: 2c_n) \in \cL$ if and only if $H = y_1$. We arrive at a contradiction $0 = y_1[-\lambda^2 x_0 + x_1 + \lambda x_2] \neq 0$. 

Second, it remains to prove that $\rk(F_i) \geq 4(n-1)$ for $i \in\{ 5,9\}$. We consider $F_5$, and the arguments are analogous to $F_9$. We argue by contradiction: assume $Z \subset \PP^n$ is a set of $|Z| \leq 4(n-1)-1$ such that $I(Z) \subset \Ann(F_5)$ (Lemma \ref{lemma: apolarity}). By Corollary \ref{Corollary: geometric conse generalized lower bound}(i), we can assure that $Z = Z_1 \cup Z_2 \cup Z_3 \cup Z_4$ with $Z_i$ a set of $n-1$  linearly independent  points such that $\langle Z_{i} \rangle \cap \cL = \emptyset$ for every $i \in\{ 1,2,3\}$. Hence, $|Z_4| = n-2$ and we can take a hyperplane $\Var(H_c)$ containing $\langle Z_4 \rangle + \cL$. We will show that for any such hyperplane $\Var(H_c)$, the point $Q= H_{c}[q]^{\vee} = (0,c_2,0, c_3, \hdots, c_n) \notin \cL$. Therefore, $I(Z) \not\subset \Ann(F)$ by Corollary \ref{Corollary: geometric conse generalized lower bound}(ii), and we arrive at a contradiction. If $Q \in \cL$, then $H_{c} = y_2$, and we get $H_{c}[F] = x_0x_1^2$. We define hyperplanes $\Var(H_{1}) = \langle Z_1 \cup \{x_1^{\vee}\} \rangle$, $\Var(H_{2}) = \langle Z_2 \cup \{x_0^{\vee}\} \rangle$ and $\Var(H_{3}) = \langle Z_3 \cup \{x_0^{\vee}\} \rangle$. By construction, $H_{3}H_{2}H_{1}H_{c} \in I(Z) \subset \Ann(F)$ but $H_{3}H_{2}H_{1}H_{c}[F] = H_{3}H_{2}H_{1}[x_0x_1^2] \neq 0$. The proof of the proposition is complete. 
\end{proof}

We end this subsection with two observations. 

\begin{remark} (i) Studying closely $\Ann(F_i)_2$ for $i \in\{ 4,\hdots,10\}$, one can see that $\dim\Ann(F_i)_2 \leq \binom{n}{2}+1$. Thus, 
Remark \ref{remark: catalecticant bound} only gives $\rk(F_{i}) \geq 2n$, while our new lower bound is $\rk(F_i) \geq 4(n-1)$. 

\vspace{0.15cm} 
\noindent (ii) Theorem \ref{Proposition: generalized lower bound} and Corollary \ref{Corollary: geometric conse generalized lower bound} impose constrictions for an ideal $I(Z)$ of $|Z| = 4(n-1)$ points to be contained in $\Ann(F_i)$. First, $Z = Z_1 \cup Z_2 \cup Z_3 \cup Z_4$ where each $Z_j$ is a set of $n-1$  linearly independent  points such that $\langle Z_j \rangle \cap \cL_i = \emptyset$ for every $j \in\{ 1,2,3\}$. Then, also $Z_4$ is a set of $n-1$  linearly independent  points such that $\langle Z_4 \rangle \cap \cL_i = \emptyset$. Such set splittings can always be constructed iteratively for each choice of $Z_1$. 
\end{remark}

\subsection{Upper bounds for the symmetric rank}
\label{Section: Upper bound}

In this subsection, we provide upper bounds for the quartics $F_i$. Consequently, we establish that $4(n-1) \leq \rk(F_i) \leq 4(n-1)+2$ for every $i \in\{4,\hdots,10\}$ (Theorem \ref{Prop: upper bounds for reducible quartics}). In particular, we prove that the lower bound $4(n-1)$ is optimal for some of these forms, and we will observe that it may not be the case for the others, at least for low values of $n$. 

\vspace{0.15cm}
We start showing that $\rk(x_0x_1(x_2^2 + \cdots + x_n^2)) = \rk(x_0x_1(x_0^2 + \cdots + x_n^2) = 4(n-1)$, key points in the main result of this subsection. 
\begin{proposition} \label{Prop: upper bound for F4} If $n \geq 3$, then 
$\rk(x_0x_1(x_2^2+\cdots+x_n^2)) = 4(n-1)$. Moreover, it has infinitely many decompositions of length $4(n-1)$. 
\end{proposition}

\begin{proof} We denote $F = x_0x_1(x_2^2+\cdots+x_n^2)$. Thanks to Proposition \ref{Theorem: lower rank of quartics}, we only have to prove the upper bound. We will construct a set $Z\subset \PP^n$ of $4(n-1)$ points such that $I(Z) \subset \Ann(F)$. The idea is to {\em lift and glue} decompositions of the cubics $x_0(x_2^2 + \cdots + x_n^2)$ and $x_1(x_2^2 + \cdots + x_n^2)$. 

An explicit decomposition of length $2(n-1)$ of the cubic $x_1(x_2^2 + \cdots + x_n^2)$ is given in \cite[Proposition 7.2(1)]{Landsberg-Teitler}. Thus, it defines an ideal $I(Z_1) \subset \PP^{n-1}$ of $2(n-1)$ points with $I(Z_1) \subset \Ann(x_1(x_2^2 + \cdots + x_n^2)) \subset S$. It is straightforward to check that $Z_1 \cap \Var(y_1) = \emptyset$. For each $p = (a_1: a_2: \cdots : a_n) \in Z_1$, we define $p_{y_0-y_1} = (a_1: a_1: a_2: \cdots : a_n)$ and let $Z_{1,y_0-y_1} \subset \PP^n$ be the set of such $2(n-1)$ points. Analogously, we define $p_{y_0+y_1} = (-a_1: a_1: a_2: \cdots : a_n)$ and $Z_{1,y_0+y_1}$ the associated set of $2(n-1)$ points. We set  $Z = Z_{1,y_0-y_1} \cup Z_{1,y_0+y_1} \subset \PP^n$. Since $Z_1 \cap \Var(y_1) = \emptyset$, the set $Z$ is indeed a reduced set of $4(n-1)$ points in $\PP^n$. In particular, $Z \subset \Var(y_0-y_1) \cup \Var(y_0+y_1)$. 

We will show $I(Z) \subset \Ann(F)$. Since $y_0^2 - y_1^2 \in \Ann(F)$ and
\[\Ann(F) = (y_0^2) + (\Ann(y_1(y_2^2+ \cdots + y_n^2)) \cap S,\]
it is enough to show that  $I(Z) = (y_0^2 - y_1^2) + (I(Z_1)) \subset \Ann(F)$, where $I(Z_1) \subset S$ denotes the ideal of $Z_1 \subset \PP^{n-1}$ and $(I(Z_1))$ denotes the ideal generated by $I(Z_1)$ in $S$. 

First, we notice that $\Var(I(Z)) = Z$, so it only remains to prove that $I(Z)$ is a radical ideal. To do it, we  will see that 
 $I(Z_1)_{y_0-y_1} = (y_0 - y_1) + (I(Z_1))$ (respectively 
$I(Z_1)_{y_0+y_1} = (y_0 + y_1) + (I(Z_1))$) and $I(Z) = I(Z_1)_{y_0-y_1} \cap I(Z_1)_{y_0+y_1}.$

Observe that $\Var((y_0 - y_1) + (I(Z_1))) = Z_{1,y_0-y_1}$ and let $g \in I(Z_1)_{y_0-y_0}$. We write
\[g = g(y_0 - y_1 + y_1, y_1, y_2, \hdots, y_n) = h(y_0, \hdots, y_n)(y_0-y_1) + f(y_1, \hdots, y_n).\] 
 Since $g$ and $y_0-y_1$ belong to $I(Z_1)_{y_0-y_1}$ and $f \in S$, we have that $f$ vanishes on $Z_{1} \subset \PP^{n-1}$ and we get $f \in I(Z_1)$. An analogous argument proves that $I(Z_1)_{y_0+y_1} = (y_0 + y_1) + (I(Z_1))$. 

Now, assume $g \in I(Z_1)_{y_0-y_1} \cap I(Z_1)_{y_0+y_1}$ and, as before, we write
\[g = h(y_0,\hdots,y_n)(y_0-y_1) + f(y_1,\hdots,y_n) = h'(y_0,\hdots,y_n)(y_0+y_1) + f'(y_1,\hdots,y_n),\]
giving $f,f' \in I(Z_1)$. We have $h(y_0-y_1) =  h'(y_0 + y_1) + f' - f \in I(Z_1)_{y_0+y_1}$. By construction, $\Var(y_0 - y_1) \cap Z_{1,y_0+y_1} = \emptyset$, so $h \in I(Z_1)_{y_0+y_1}$. Writing $h = h''(y_0+y_1) + f'' \in (y_0+y_1) + (I(Z_1))$, we get 
\[g = (y_0^2 - y_1^2)h''h + f''(y_0-y_1) + f \in I(Z) = (y_0^2 - y_1^2) + (I(Z_1)),\]
and the proof is complete.  

Moreover, $\Ann(F)$ contains any form $y_0^2 - \lambda^2 y_1^2 = (y_0 + \lambda y_1)(y_0 - \lambda y_1)$, $\, \lambda \in \CC$. Analogously as above, for any $\lambda \neq 0$, we can construct different sets $Z_{\lambda} = Z_{1,y_0-\lambda y_1} \cup Z_{1,y_0+\lambda y_1} \subset \Var(y_0-\lambda y_1) \cup \Var(y_0+\lambda y_1)$. Similar arguments shows that $I(Z_\lambda) \subset \Ann(F)$. 
\end{proof}

\begin{proposition}\label{propos: upper bound for F80}
If $n \geq 2$, then $\rk(x_0x_1(x_0^2 + x_1^2 + \cdots + x_n^2)) = 4(n-1)$. It has infinitely many decompositions of length $4(n-1)$. 
\end{proposition}
\begin{proof}
Let $F = x_0x_1(x_0^2 + x_1^2 + x_2^2 + \cdots + x_n^2)$. For $n=2$, we prove the bound by explicitly producing the ideal of four points
\[(y_0^2 - \lambda^2 y_1^2, y_1^2 - 3y_2^2) \subset \Ann(x_0x_1(x_0^2 + x_1^2 + x_2^2)), \; 0 \neq \lambda \in \CC;\]
and for $n=3$, the ideal of eight points 
\[(y_0^2 - \lambda^2 y_1^2, y_1^2 - 3y_2^2, y_2^2 - y_3^2) \subset \Ann(x_0x_1(x_0^2 + x_1^2 + x_2^2 + x_3^2)), \; 0 \neq \lambda \in \CC.\]
Now assume that $n \geq 4$ and write 
$F = x_0x_1(x_0^2 + x_1^2 + x_2^2) + x_0x_1(x_3^2 + \cdots + x_{n}^2)$. By Proposition \ref{Prop: upper bound for F4}, we have $\rk(x_0x_1(x_3^2 + \cdots + x_n^2)) \leq  4(n-2)$. 
Thus, for $n \geq 4$, we obtain $\rk(F_{8,0}) \leq 4 + 4(n-2) = 4(n-1)$, and the result follows.
\end{proof}

\begin{lemma}\label{Lemma: Keples ternary quartics} Let $G \in \CC[x_0, x_1, x_2]_4$ be a concise form such that $\dim \Ann(G)_2 \leq 2$. 
\begin{itemize}
       \item[(i)] If $\Var(\Ann(G)_2) \subset \PP^2$ is a finite set of $k$ points, then $\rk(G) = 4$ if $k = 4$ and $\rk(G) = 6$ otherwise.  
       \item[(ii)] If $\Ann(G)_2 = (0)$, then $\rk(G) = 6$.
\end{itemize}
\end{lemma}
\begin{proof} 
The first statement follows by   \cite[Theorem 3.2]{kleppe},  while the second one is a consequence of \cite[Theorem 3.7]{kleppe}. 
\end{proof}

We are ready to present the main result of this section. 

\begin{theorem}\label{Prop: upper bounds for reducible quartics}\label{thm: rank of the 10 families} If $n \geq 3$, then
    \begin{enumerate}
        \item $\rk(F_1) = \rk(x_0^{2}(x_0^2+\cdots+x_n^2)) = 3n$.
        \item $\rk(F_2) = \rk(x_0^{2}(x_1^2+\cdots+x_n^2)) = 3n$.
        \item $\rk(F_3) = \rk(x_0^{2}(x_0x_1+x_2^2+\cdots+x_n^2)) = 3n+1$.
        \item $\rk(F_4) = \rk(x_0x_1(x_2^2+\cdots+x_n^2)) =  4(n-1)$.
        \item $4(n-1) \leq \rk(F_5) = \rk(x_0x_1(x_1x_2+x_3^2+\cdots+x_n^2)) \leq 4(n-1)+2$.
        \item $\rk(F_6) = \rk(x_0x_1(x_1^2+\cdots+x_n^2)) = 4(n-1)$.
        \item $4(n-1) \leq \rk(F_{7}) = \rk(x_0(x_0+x_1)(x_1^2 + \cdots + x_n^2)) \leq 4(n-1)+2$.
        \item $4(n-1) \leq \rk(F_{8,\lambda}) = x_0(\lambda x_0+x_1)(x_0^2+\cdots+x_n^2)) \leq 4(n-1)+2$ for $0 \neq \lambda \in \CC$ and $\rk(F_{8,0}) = \rk(x_0x_1(x_0^2 + x_1^2 + x_2^2 + \cdots + x_n^2)) = 4(n-1)$.
        \item $4(n-1) \leq \rk(F_9) = \rk(x_0x_1(x_0^2 + x_1x_2 + x_3^2 + \cdots + x_n^2)) \leq 4(n-1)+2$.
        \item $4(n-1) \leq \rk(F_{10,\lambda})=\rk(x_0(-\lambda^2 x_0+x_1+\lambda x_2)(x_0x_1+x_2^2+\cdots+x_n^2)) \leq 4(n-1)+2$ for every $\lambda \in \CC$. 
    \end{enumerate}
    Moreover, $F_4, F_6$ and $F_{8,0}$ admits infinitely many decompositions of length $4(n-1)$. 
\end{theorem}
\begin{proof} The lower bounds are proved in Proposition \ref{Theorem: lower rank of quartics}. The ranks of $F_1, F_2, F_3$, are given in Proposition \ref{Proposition: rank of high multiplicity hyperplane}, and the ranks of $F_4$ and $F_{8,0}$ are given in Propositions \ref{Prop: upper bound for F4} and \ref{propos: upper bound for F80}. It remains to prove the upper bounds of the remaining quartic forms. We proceed by expressing them as a sum $G_i + H_i$, with $G_i \in \CC[x_0,x_1,x_2]_4$ as in Lemma \ref{Lemma: Keples ternary quartics} and $H_i$ is of type $F_{8,0}$ as in Proposition \ref{propos: upper bound for F80}. 
\begin{itemize}
\item[(5)] Write $F_5 = x_0x_1(-x_0^2 - x_1^2 + x_1x_2) + x_0x_1(x_0^2 + x_1^2 + x_3^2 + \cdots + x_n^2) = G_5 + H_5$. The rank of $G_5$ is $6$ by Lemma \ref{Lemma: Keples ternary quartics}, indeed $\Ann(G_5)_2 = (y_0^2 + 3y_1y_2, y_2^2)$ and $\Var(\Ann(G_5)_2) = \{(0:1:0)\}$. By Proposition \ref{propos: upper bound for F80}, we have $\rk(H_5) = 4(n-2)$. Hence, 
$\rk(F_5) \leq 6 + 4(n-2) = 4(n-1)+2$. 

\item[(6)] Write $F_6 = x_0x_1(-x_0^2-x_1^2+x_2^2) + x_0x_1(x_0^2+2x_1^2 + x_3^2 + \cdots + x_n^2) = G_6 + H_6.$
The ideal $(y_1^2 + 3y_2^2, y_0^2 + 3y_2^2) \subset \Ann(G_6)$ is the radical ideal of $4$ points, so $\rk(G_6) = 4$ by Lemma \ref{Lemma: Keples ternary quartics}. On the other hand, $\rk(H_6) \leq 4(n-2)$ by Proposition \ref{propos: upper bound for F80}, hence $\rk(F_6) \leq 4 + 4(n-2) = 4(n-1).$
\item[(7)] Write $F_7 = x_0(x_0+x_1)(-2x_0^2 - 2x_0x_1 + x_2^2) + x_0(x_0+x_1)(x_0^2 + (x_0+x_1)^2 + x_3^2 + \cdots + x_n^2) = G_7 + H_7$. Since $\Ann(G_7)_2 = (0)$,  $\rk(G_7) = 6$ by Lemma \ref{Lemma: Keples ternary quartics}, while $\rk(H_7) =4(n-2)$ by Proposition \ref{propos: upper bound for F80}, hence $\rk(F_7) \leq 6 + 4(n-2) = 4(n-1)+2$. \end{itemize} 
We conclude cases (8) to (10) with similar arguments as in (7), writing:
\[\begin{array}{lllllllll}
F_{8,\lambda} \!  & = & x_0(\lambda x_0+x_1)(-\lambda^2x_0^2 - 2\lambda x_0x_1 + x_2^2) + \\
& & x_0(\lambda x_0+x_1)(x_0^2 + (\lambda x_0 + x_1)^2 + x_3^2 + \cdots + x_n^2)  = G_{8,\lambda} + H_{8,\lambda},
\end{array}\] 
with $0 \neq \lambda \in \CC$;
\[\begin{array}{llllllll}
  F_9 & = & x_0x_1(-x_1^2 + x_1x_2) + x_0x_1(x_0^2 + x_1^2 + x_3^2 + \cdots + x_n^2) = G_9 + H_9, \\
\end{array}\]
and
\[\begin{array}{llllll}
F_{10,\lambda} & = & x_0(-\lambda^2 x_0 + x_1 + \lambda x_2)(-x_0^2 - (-\lambda^2 x_0 + x_1 + \lambda x_2)^2 + x_0x_1 + x_2^2)] + \\
& & x_0(-\lambda^2 x_0 + x_1 + \lambda x_2)(x_0^2 + (-\lambda^2 x_0 + x_1 + \lambda x_2)^2 + x_3^2 + \cdots + x_n^2) \\
& = & G_{10,\lambda} + H_{10,\lambda}    
\end{array}\]
with $\lambda \in \CC$. 
\end{proof}

 Next, we discuss further the symmetric rank of $F \in \{F_7,F_{8,\lambda \neq 0}, F_9, F_{10,\lambda}\}$, $\lambda \in \CC$. 
For $n = 3$, Theorem \ref{Prop: upper bounds for reducible quartics} assures $8 \leq \rk(F) \leq 10$. Actually, we have $\rk(F) > 8$. Indeed, one can check that 
\[\Ann(F)_2 = \langle y_2y_3, y_2^2 - y_3^2 \rangle,\]
which contains a square $(y_2 + i y_3)^2$, where $i^2 = -1$. If $Z \subset \PP^3$ is a set of $8$ points such that $I(Z) \subset \Ann(F)$, then $\dim I(Z)_2 \geq 2$ and we obtain $I(Z)_2 = \Ann(F)_2$. Since $I(Z)$ is a radical ideal, $y_2 + iy_3 \in I(Z) \subset \Ann(F)$, and this contradicts the conciseness of $F$. The picture changes for $n \geq 4$. We still have 
\[\Ann(F)_2 = \Ann(x_2^2 + \cdots + x_n^2)_2
\subset \CC[y_2,\hdots, y_n].\] 
Notwithstanding, $\Ann(F)_2$ is now big enough so that we cannot rule out that it may contain the ideal of $4(n-1)$ or $4(n-1)+1$ points. Thus, the behavior of $\rk(F)$ may differ for low values of $n$ and asymptotically.

\begin{remark}\label{remark: identifiability}
By \cite{genericforms}, the rank of a general quartic in $n+1$ variables is 10 for $n=3$, it is 15 for $n=4$ and it is $\ru{\frac{1}{n+1}\binom{n+4}{4}}$ for every $n\ge 5$, so Theorem \ref{thm: rank of the 10 families} implies that our families have symmetric rank smaller than the generic rank, except for a few small values of $n$ (compare Remark \ref{remark: examples of supergeneric and maximal rank} for a discussion of these small cases). However, our reducible quartics exhibit another unexpected behaviour
. Indeed, \cite[Theorem 1.1]{COVsubgeneric} states that for $n=3$, a general rank-8 quaternary quartic has exactly two Waring decompositions of length 8, while for $n\ge 4$, a general quartic of subgeneric rank has a unique rank decomposition. On the contrary, $F_4$, $F_6$ and $F_{8,0}$ have infinitely many rank decompositions.
\end{remark}

A final note regarding the partially symmetric rank. If $i \in\{ 4,\hdots,10\}$, then we have established that $\prk(F_i) \geq 3(n-1)$ and $\rk(F_i) \geq 4(n-1)$, with equality for some of them (Propositions \ref{Theorem: lower rank of quartics}, \ref{Prop: upper bound for F4} and \ref{propos: upper bound for F80}). The core of the argument for the symmetric rank is that if $Z \subset \PP^n$ is an ideal of $4(n-1)-1$ points with $I(Z) \subset \Ann(F_i)$, then by distributing the points on hyperplanes, we construct a form $H_1\cdots H_4 \in I(Z)$ which cannot belong to $\Ann(F)_4$, hence a contradiction. However, this does not yield a contradiction if assuming $I(Z) \subset \Ann(J_{F_i})$ given that $\Ann(J_{F_i})_4 = S_4$. Nevertheless, it would be very interesting to know if $\prk(F_i) = \rk(F_i)$ and if Comon's conjecture holds for this type of quartics.
Motivated by these observations, we end this section by posting the following problems and questions.

\begin{enumerate}
    \item For $n = 3$ and $F \in \{F_7, F_{8,\lambda \neq 0}, F_9, F_{10,\lambda}\}$, we know that $9 \leq \rk(F) \leq 10$. If $\rk(F) = 9$, then we would obtain $\rk(F) \geq 4(n-1)+1$ for $n \geq 3$. Which is the precise value of $\rk(F)$? Is $\rk(F) > 4(n-1)$ for $n \geq 4$?
    \item Is $\rk(F_5) = 8$ for $n = 3$? This would imply that $\rk(F_5) = 4(n-1)$ for $n \geq 3$. 
    \item Is $\rk(F_i) = 4(n-1)$ asymptotically in $n$?
    \item Do we have an equality $\prk(F_i) = \rk(F_i)$ for $i\in\{ 4,\hdots,10\}$? As pointed out above, our arguments do not discourage a strict inequality. 
\end{enumerate}

\bibliographystyle{plain}
\bibliography{refcomonimpan.bib}

\end{document}